\documentclass[11pt]{article}
\usepackage[a4paper]{anysize}\marginsize{2.0cm}{2.0cm}{1.0cm}{1.0cm}
\pdfpagewidth=\paperwidth \pdfpageheight=\paperheight
\usepackage{pgf,tikz,hyperref}
\usetikzlibrary{arrows}
\usepackage{amsfonts,amssymb,amsthm,amsmath,eucal}
\usepackage{graphicx}
\usepackage{ltablex}
\usepackage{caption, booktabs,verbatim}
\theoremstyle{plain}
\newtheorem{thm}{Theorem}[section]
\newtheorem{theorem}[thm]{Theorem}

\newtheorem{lemma}[thm]{Lemma}
\newtheorem{proposition}[thm]{Proposition}
\newtheorem{corollary}[thm]{Corollary}

\theoremstyle{definition}
\newtheorem{definition}[thm]{Definition}
\newtheorem{remark}[thm]{Remark}
\newtheorem{example}[thm]{Example}

\newtheorem{thevarthm}[thm]{\varthmname}

\newenvironment{varthm*}[1]{\trivlist\item[]{\bf #1.}\it}{\endtrivlist}

\def\cftil#1{\ifmmode\setbox7\hbox{\accent"5E#1}\else\setbox7\hbox{\accent"5E#1}\penalty 10000\relax\fi\raise 1\ht7\hbox{\lower1.1ex\hbox to 1\wd7{\hss\accent"7E\hss}}\penalty 10000\hskip-1\wd7\penalty 10000\box7 }
\newcommand\eps{\varepsilon}
\newcommand\be{\begin{eqnarray*}}
\newcommand\ee{\end{eqnarray*}}

\newcommand\C{\mathbb C}
\newcommand\call{\mathcal L}
\newcommand\calm{\mathcal M}

\renewcommand\P{\mathbb P}

\newcommand\newop[2]{\def#1{\mathop{\rm #2}\nolimits}}
\newop\edim{edim}
\newop\Zeroes{Zeroes}
\newop\Jac{Jac}
\newop\Ass{Ass}
\newop\SL{SL}
\newop\PGL{{\P}GL}
\newop\Km{Km}
\newop\reg{reg}
\newop\Hess{Hess}
\newop\tr{tr}
\newop\Id{Id}

\newcommand\keywords[1]{{\renewcommand\thefootnote{}\footnotetext{\textit{Keywords:} #1.}}}
\newcommand\subclass[1]{{\renewcommand\thefootnote{}\footnotetext{\textit{Mathematics Subject Classification (2020):} #1.}}}

\usepackage{authblk}
\begin{document}
\definecolor{uuuuuu}{rgb}{0.26666666666666666,0.26666666666666666,0.26666666666666666}

\author{Luca Chiantini, \L ucja Farnik, Giuseppe Favacchio, Brian Harbourne, Juan Migliore,\\ Tomasz Szemberg, Justyna Szpond}

\title{On the classification of certain geproci sets}
\date{\today}
\maketitle
\thispagestyle{empty}

\begin{abstract}
In this short note we develop new methods toward the ultimate goal of classifying geproci sets in $\mathbb P^3$. We apply these methods to show that among sets of $16$ points distributed evenly on $4$ skew lines, up to projective equivalence there are only two distinct geproci sets. We give different geometric distinctions between these sets. The methods we develop here can be applied in a more general set-up; this is the context of  follow-up work \cite{NEW}.
\end{abstract}

\subclass{MSC 14M05 \and MSC 14M10 \and MSC 14N05 \and MSC 14N20}

%*****************************************************************************

\section{Introduction}

The study of geproci sets was initialized in one of the previous workshops on Lefschetz properties held in Levico Terme (Italy) in 2018. In the present note we work exclusively over the field $\C$ of complex numbers.
\begin{definition}[A geproci set of points]
   We say that a set of points $Z\subset\P^N_{\mathbb{C}}$ with $N\geq 3$ is \emph{geproci} (for GEneral PROjection is a Complete Intersection), if its general projection to a hyperplane is a complete intersection.
\end{definition}
We say that a geproci set is \emph{trivial}, if it is already contained in a hyperplane. From now on we consider only nontrivial geproci sets.

So far nontrivial geproci sets have been discovered only in $\P^3$. They project to a plane, where their images are the intersection points of two curves of degrees $a$ and $b$ (equivalently: the ideal of the projection has exactly two generators: one of degree $a$ and another of degree $b$). We assume that $a\leq b$ and we refer to such sets of points as $(a,b)$-geproci.
\begin{example}[A grid]\label{exp: grid}
   Assume that we have two positive integers $a\leq b$. Let $Z\subset\P^3$ be a grid, i.e., the set of all intersection points among lines in two sets $\call=\left\{L_1,\ldots,L_a\right\}$ and $\calm=\left\{M_1,\ldots,M_b\right\}$ such that lines from the same set, either $\call$ or $\calm$, are pairwise skew but any two lines from distinct sets intersect in a point. It is elementary to see that a grid is an $(a,b)$-geproci set for all values of $a,b$ and the set is nontrivial for $a,b\geq 2$ and it is contained in a unique quadric for $a,b\geq 3$.
\end{example}
Grids exist for any values of $a$ and $b$. They were studied extensively in \cite{ChiMig21}. Here we are interested in geproci sets which are not grids but half grids.
\begin{definition}[A half grid]
    We say that a nontrivial $(a,b)$-geproci set $Z\subset\P^3$ is a \emph{half grid}, if it is not a grid but one of the curves determining its general projection as a complete intersection can be taken as a union of lines.
\end{definition}
\begin{remark}
   General projections of points can be collinear only if the points are collinear before. So an $(a,b)$ half grid in $\P^3$ consists either of $a$-tuples of $b$ points each distributed on $a$ lines or of $b$-tuples of $a$ points distributed on $b$ lines. A priori it could be that the lines containing points of $Z$ intersect. However \cite[Proposition 4.14]{BIG} provides an easy argument that this is not the case. Indeed, removing all lines but three from the half grid, we obtain a $(3,a)$ (or a $(3,b)$) geproci set. This set must be a grid, so that in particular the lines must be disjoint.
\end{remark}
In particular half grids are geproci by definition. Most of the known geproci sets which are not grids are half grids. As a matter of fact, we know at the time of this writing only three exceptions -- see the introduction to \cite{BIG} for details. This justifies that our interest here focuses on half grids.

Our main result in this note is the following.
\begin{theorem}[Classification of $(4,4)$ half grids]\label{thm: main}
   Let $Z\subset\P^3$ be a $(4,4)$ half grid. Then, up to projective change of coordinates, $Z$ is either
   \begin{itemize}
       \item[A)] the anharmonic case (see Section \ref{ssec:anharmonic})
$$\begin{array}{llllllll}
(1:0:0:0), && (0:1:0:0), && (1:1:0:0), &&  (0:1:1:1),\\
(0:0:1:0), && (0:0:0:1), && (0:0:1:1), && (1:1:1:0),\\
(1:0:1:0), && (0:1:0:1), && (1:1:1:1), && (1:0:0:-1),\\
(1:0:\eps:0), && (0:1:0:\eps), && (1:1:\eps:\eps), && (1: \eps: \eps: \eps-1),
\end{array}$$
where $\eps$ is a primitive root of unity of order six, or
       \item[B)] the harmonic case (see Section \ref{ssec:harmonic})
$$\begin{array}{llllllll}
(1:0:0:0), && (0:1:0:0), && (1:1:0:0), &&  (1:0:0:-1),\\
(0:0:1:0), && (0:0:0:1), && (0:0:1:1), && (0:1:1:0),\\
(1:0:1:0), && (0:1:0:1), && (1:1:1:1), && (1:1:1:-1),\\
(1:0:-1:0), && (0:1:0:-1), && (1:1:-1:-1), && (-1:1:1:1).
\end{array}$$
   \end{itemize}
\end{theorem}
The points in the statement of Theorem \ref{thm: main} are organized so that the half grid property is immediately visible: points in the columns are collinear. The main difference between both cases is that in the harmonic case there are exactly $4$ lines containing $4$ of configuration points, whereas in the anharmonic case additional collinearity can be observed for the $4$ points in the bottom row.

In order to put this result in some perspective, let us recall that in \cite[Theorem 4.10]{BIG} we proved that the only non-grid $(3,4)$-geproci set is the half grid determined by points in the $D_4$ root system. Up to projective change of coordinates, its $12$ points can be listed explicitly as
$$\begin{array}{llllll}
(1:1:0:0),& (1:0:1:0),& (1:0:0:1),& (0:1:1:0),& (0:1:0:1),& (0:0:1:1),\\
(1:-1:0:0),& (1:0:-1:0),& (1:0:0:-1),& (0:1:-1:0),& (0:1:0:-1),& (0:0:1:-1).
\end{array}$$
See Theorem \ref{thm: (3,4)-geproci} for a more precise statement.
%*****************************************************************************
\section{Inputs from projective geometry}
%*****************************************************************************
To fix notation let $(x:y:z:w)$ be projective coordinates on $\P^3$.
We denote by $j(P_1,P_2;P_3,P_4)$ the cross-ratio of an ordered set of four collinear points (see \cite{BIG} for the definition and first results). For integers $\left\{i,j,k,l\right\}=\left\{1,2,3,4\right\}$, we write $(i,j,k,l)$ to indicate the permutation, which sends $1$ to $i$, $2$ to $j$, $3$ to $k$ and $4$ to $l$. Note that this is not cycle notation!

 It is well-known that the cross-ratio is invariant under the Klein group, namely
$$j(P_1,P_2;P_3,P_4)=j(P_{\sigma(1)},P_{\sigma(2)};P_{\sigma(3)},P_{\sigma(4)}),$$
where $\sigma$ is one of the following permutations:
\begin{equation}\label{eq: standard c-r}
    (1,2,3,4),\; (2,1,4,3),\; (3,4,1,2),\; (4,3,2,1).
\end{equation}
Note that all elements in the Klein group are involutions.

There are two exceptional cases under which there are additional permutations leaving the cross-ratio invariant:\\[1em]
\textbf{The harmonic case.} If $j(P_1,P_2;P_3,P_4)\in\left\{-1, 1/2, 2\right\}$, then the following permutations leave the cross-ratio invariant:
\begin{equation}\label{eq: harmonic}
\begin{array}{llll}
(1,2,3,4),& (2,1,4,3),& (3,4,1,2),& (4,3,2,1),\\
(1,2,4,3),& (2,1,3,4),& (3,4,2,1),& (4,3,1,2).
\end{array}\end{equation}
\noindent
\textbf{The anharmonic case.} If $j(P_1,P_2;P_3,P_4)\in\left\{\frac{1}{2}+\frac{\sqrt{3}}{2}i, \frac{1}{2}-\frac{\sqrt{3}}{2}i\right\}$, i.e., it is a primitive root of unity of order $6$, then the following permutations leave the cross-ratio invariant:
\begin{equation}\label{eq: anharmonic}
\begin{array}{llll}
(1,2,3,4),& (2,1,4,3),& (3,4,1,2),& (4,3,2,1),\\
(1,3,4,2),& (2,4,3,1),& (3,1,2,4),& (4,2,1,3),\\
(1,4,2,3),& (2,3,1,4),& (3,2,4,1),& (4,1,3,2).
\end{array}\end{equation}

It is a classical fact in projective geometry (see, e.g., \cite[Paragraph 3.4.1]{EisHar16}) that given four lines in $\P^3$  not on a quadric surface, there are two (counted with multiplicities) transversals to these lines.

We saw in Example \ref{exp: grid} that any $(a,b)$ grid with $a,b\geq 3$ is contained in a quadric.
Our first result here is a criterion when lines determining a $(2,4)$ grid are contained in a quadric.
\begin{lemma}[Quadrics and $(2,4)$ grids.]\label{lem: quadrics and (2,4) grids}
   Let $R, R'\subset\P^3$ be a pair of skew lines. Let $P_1,\ldots,P_4$ be a set of mutually distinct points on $R$ and let $P_1',\ldots,P_4'$ be a set of mutually distinct points on $R'$. Let $r_i$ be the line determined by $P_iP_i'$ for $i=1,\ldots,4$.
   The lines $r_1,\ldots,r_4$ are contained in a quadric if and only if
   $$j(P_1,P_2;P_3,P_4)=j(P_1',P_2';P_3',P_4').$$
\end{lemma}
\begin{proof}
The lines $r_1,r_2,r_3$ are pairwise skew, so they determine a unique quadric $Q$. This quadric contains $R$ and $R'$ because it has at least $3$ points common with both lines.
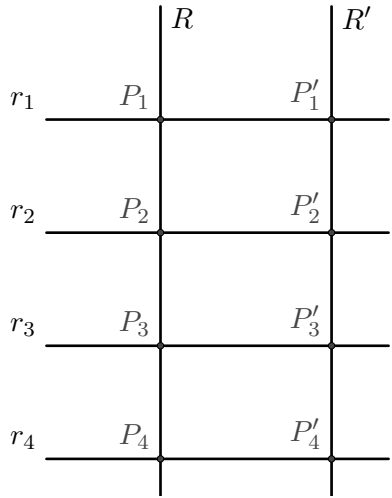
\begin{figure}[h!]
    \centering
\begin{tikzpicture}[line cap=round,line join=round,>=triangle 45,x=1cm,y=1cm,scale=0.15]
\clip(-25.8109674696541,-13.19644082133315) rectangle (41.76413067908503,34.5997794785621);

%%%%%%%%%%%%%%%%%%%%%%%%%%%%%%% PROSTE POZIOME %%%%%%%%%%%%%%%%%%%%%%%%%%%%%%%%%%%%%
%%% r_1
\draw [line width=1pt] (-10,20)-- (20,20);
\draw (-10,20) node[anchor=south east] {$r_1$};
%%% r_2
\draw [line width=1pt] (-10,10)-- (20,10);
\draw (-10,10) node[anchor=south east] {$r_2$};
%%% r_3
\draw [line width=1pt] (-10,0)-- (20,0);
\draw (-10,0) node[anchor=south east] {$r_3$};
%%% r_4
\draw [line width=1pt] (-10,-10)-- (20,-10);
\draw (-10,-10) node[anchor=south east] {$r_4$};

%%%%%%%%%%%%%%%%%%%%%% PROSTE PIONOWE %%%%%%%%%%%%%%%%%%%%%%%%%%%%%%%%%%%%%%%%%%%%%%

%%% R_b
\draw [line width=1pt] (0,30)-- (0,-40);
\draw (0,29) node[anchor= west] {$R$};
%%% R_c
\draw [line width=1pt] (15,30)-- (15,-40);
\draw (15,29) node[anchor= west] {$R'$};

%%%%%%%%%%%%%%%%%%%%%%%%%%%%%%%  PUNKTY %%%%%%%%%%%%%%%%%%%%%%%%%%%%%%%%%%%%
\draw [fill=uuuuuu] (0,20) circle (8pt);
\draw[color=uuuuuu] (0,20) node[anchor=south east] {$P_1$};

\draw [fill=uuuuuu] (0,10) circle (8pt);
\draw[color=uuuuuu] (0,10) node[anchor=south east] {$P_2$};

\draw [fill=uuuuuu] (0,0) circle (8pt);
\draw[color=uuuuuu] (0,0) node[anchor=south east] {$P_3$};

\draw [fill=uuuuuu] (0,-10) circle (8pt);
\draw[color=uuuuuu] (0,-10) node[anchor=south east] {$P_4$};

\draw [fill=uuuuuu] (15,20) circle (8pt);
\draw[color=uuuuuu] (15,20) node[anchor=south east] {$P_1'$};

\draw [fill=uuuuuu] (15,10) circle (8pt);
\draw[color=uuuuuu] (15,10) node[anchor=south east] {$P_2'$};

\draw [fill=uuuuuu] (15,0) circle (8pt);
\draw[color=uuuuuu] (15,0) node[anchor=south east] {$P_3'$};

\draw [fill=uuuuuu] (15,-10) circle (8pt);
\draw[color=uuuuuu] (15,-10) node[anchor=south east] {$P_4'$};

\end{tikzpicture}
    \caption{A $(2,4)$ grid}
    \label{fig: (2,4) grid}
\end{figure}
If the line $r_4$ is contained in $Q$, then the two cross ratios are equal.

For the other direction, observe that there exists a unique line $r$ in the ruling of $Q$ determined by $r_1$ which passes through $P_4$. This line meets $R'$ in the unique point $P$ such that  the cross ratios $j(P_1, P_2; P_3, P_4)$ and $j(P_1', P_2'; P_3', P)$ are equal. But this implies $P=P_4'$ and we are done.
\end{proof}
We need the following simple fact about projectivities of $\P^1$. We include it here with a proof, as it was difficult to track it down in the literature.
\begin{lemma}\label{lem: projectivity of P1 with 1 FP}
    Let $\varphi$ be a projective transformation of $\P^1$ with exactly $1$ fixed point $P$. Then $\varphi$ has no other finite orbit but that of $P$.
\end{lemma}
\begin{proof}
    Without loss of generality we may assume that $P=(1:0)$. Let $M$ be a matrix representing $\varphi$. Since $\varphi(P)=P$ and $\varphi$ has no other fixed points, it has the form
    $$M=\begin{pmatrix}
    1 & \varepsilon\\
    0 & 1
    \end{pmatrix}$$
    with $\varepsilon\in\C\setminus\left\{0\right\}$, because $\varphi$ is not the identity. For a point $Q\neq P$ we have $Q=(q:1)$ for some $q\in\C$ and for a positive integer $n$ we have
    $$\varphi^n(Q)=\begin{pmatrix}
    1 & \varepsilon\\
    0 & 1
    \end{pmatrix}^n\begin{pmatrix}
        q\\ 1
    \end{pmatrix}=\begin{pmatrix}
        q+n\varepsilon\\ 1
    \end{pmatrix},$$
so the orbit of $Q$ is infinite.
\end{proof}
The next Lemma characterizes projective involutions (i.e., a projective transformations $\varphi$ with $\varphi^2=\Id$ but $\varphi\neq\Id$) with two fixed points.
\begin{lemma}[Involution of $\P^1$ with $2$ fixed points]\label{lem: involutions with 2 fixed points}
    Let $P, P'$ be two distinct points in $\P^1$. Then there exists a unique involution $\varphi$ with fixed points at $P$ and $P'$.
\end{lemma}
\begin{proof}
    We can assume that $P=(1:0)$ and $P'=(0:1)$. Then any matrix fixing (projectively) these two points has the shape
    $$M=\begin{pmatrix}
        1 & 0\\
        0 & \varepsilon
    \end{pmatrix}.$$
    The condition $M^2=\Id$ (and $M\neq\Id$) forces $\varepsilon=-1$.
\end{proof}
Our next observation is that projective transformations on two skew lines in $\P^3$ are always induced by a projective transformation of the ambient space.
\begin{lemma}[Extending projectivities on a pair of skew lines]\label{lem: projectivities on two lines}
    Let $R, R'$ be skew lines in $\P^3$. Let $\varphi$ be a projectivity on $R$ and let $\varphi'$ be a projectivity on $R'$. Then there exists a projective transformation $\Phi$ of $\P^3$, which restricts to $\varphi$ on $R$ and to $\varphi'$ on $R'$. In particular the two lines are invariant under $\Phi$.
\end{lemma}
\begin{proof}
Up to change of coordinates we can assume that $R$ is the $z=w=0$ line and $R'$ is defined by equations $x=y=0$. Let $M$ be a matrix defining $\varphi$ in $(x:y)$ coordinates and let $M'$ be a matrix defining $\varphi'$ in $(z:w)$ coordinates. Then the matrix
$$\begin{pmatrix}
    M & 0\\
    0 & M'
\end{pmatrix}$$
defines a projectivity $\Phi$ in coordinates $(x:y:z:w)$, which satisfies the requirements of the Lemma.
\end{proof}
%*****************************************************************************
\section{On the $(4,4)$ half grids}
%*****************************************************************************
In \cite[Theorem 4.10]{BIG} we showed the following classification result for $(3,4)$-geproci sets which is crucial in the sequel.
\begin{theorem}[Classification of $(3,4)$-geproci sets]\label{thm: (3,4)-geproci}
    Let $Z\subset\P^3$ be a $(3,4)$-geproci set. Then either
    \begin{itemize}
        \item[a)] $Z$ is a grid, or
        \item[b)] $Z$ is the $D_4$ configuration of points.
    \end{itemize}
\end{theorem}
Taking into account that $D_4$ does not contain any four collinear points and does not lie on a quadric surface, we derive the following immediate consequence.
\begin{corollary}[Subsets of $(4,4)$ half grids]
   Let $Z$ be a $(4,4)$ half grid and let $R_a, R_b, R_c, R_d$ be four lines covering $Z$. Then for any index $x\in\left\{a,b,c,d\right\}$ the set $Z\setminus R_x$ is a $(3,4)$ grid.
\end{corollary}
For any three mutually different symbols $x,y,z\in\left\{a,b,c,d\right\}$ we denote by $Q_{xyz}$ the quadric generated by $R_x, R_y$ and $R_z$. Since the lines are skew, these quadrics are smooth. They are also mutually distinct, because $Z$ is not contained in a quadric (it would be a grid otherwise).
\begin{figure}[h!]
    \centering
\begin{tikzpicture}[line cap=round,line join=round,>=triangle 45,x=1cm,y=1cm,scale=0.1]
\clip(-31.8109674696541,-27.19644082133315) rectangle (61.76413067908503,39.5997794785621);

%%%%%%%%%%%%%%%%%%%%%% PROSTE PIONOWE %%%%%%%%%%%%%%%%%%%%%%%%%%%%%%%%%%%%%%%%%%%%%%

%%% R_a
\draw [line width=1pt] (-13.5,30)-- (-13.5,-40);
\draw (-13.5,34) node[anchor=north] {$R_a$};
%%% R_b
\draw [line width=1pt] (0,30)-- (0,-40);
\draw (0,34) node[anchor=north] {$R_b$};
%%% R_c
\draw [line width=1pt] (15,30)-- (15,-40);
\draw (15,34) node[anchor=north] {$R_c$};

%%%%%%%%%%%%%%%%%%%%%% PROSTE UKO君E %%%%%%%%%%%%%%%%%%%%%%%%%%%%%%%%%%%%%%%%%%%%%%

%%% R_d
\draw [line width=1pt] (26,38)-- (16.7,-37);
\draw (26,40) node[anchor=north east] {$R_d$};

%%%%%%%%%%%%%%%%%%%%%%%%%%%%%%%  PUNKTY %%%%%%%%%%%%%%%%%%%%%%%%%%%%%%%%%%%%

\draw [fill=uuuuuu] (-13.5,20) circle (8pt);

\draw [fill=uuuuuu] (-13.5,10) circle (8pt);

\draw [fill=uuuuuu] (-13.5,0) circle (8pt);

\draw [fill=uuuuuu] (-13.5,-10) circle (8pt);

\draw [fill=uuuuuu] (0,20) circle (8pt);

\draw [fill=uuuuuu] (0,10) circle (8pt);

\draw [fill=uuuuuu] (0,0) circle (8pt);

\draw [fill=uuuuuu] (0,-10) circle (8pt);

\draw [fill=uuuuuu] (15,20) circle (8pt);

\draw [fill=uuuuuu] (15,10) circle (8pt);

\draw [fill=uuuuuu] (15,0) circle (8pt);

\draw [fill=uuuuuu] (15,-10) circle (8pt);

\draw [fill=uuuuuu] (24.555276497257866,26.275434462072584) circle (8pt);

\draw [fill=uuuuuu] (23.18572916907802,15.376924033599247) circle (8pt);

\draw [fill=uuuuuu] (21.816181840898172,4.478413605125911) circle (8pt);

\draw [fill=uuuuuu] (19.044287161367107,-17.87962152071977) circle (8pt);

\end{tikzpicture}
    \caption{A (4,4) half grid}
    \label{fig: (4,4) half grid}
\end{figure}
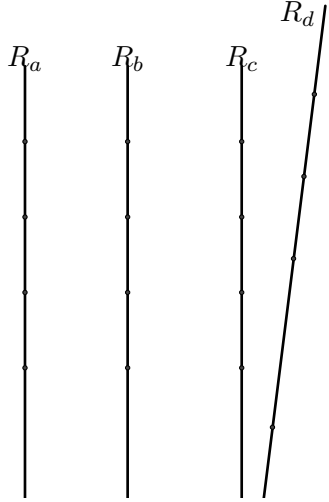

Before we impose specific coordinates on the points in $Z$, we want to determine additional collinearities. To this end we begin to label points in $Z$ on the lines provided in Figure \ref{fig: (4,4) half grid}. We begin with the line $R_c$ and we denote the points on this line with $c_1,\ldots,c_4$, see Figure \ref{fig: (4,4) half grid Rc}.
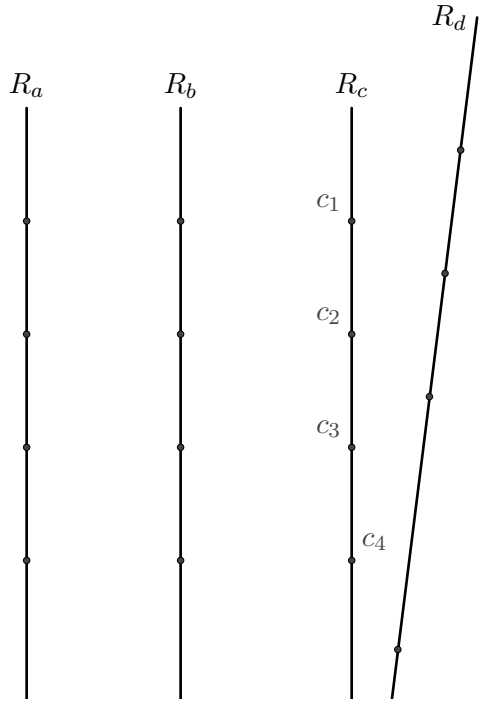
\begin{figure}[h!]
    \centering
\begin{tikzpicture}[line cap=round,line join=round,>=triangle 45,x=1cm,y=1cm,scale=0.15]
\clip(-31.8109674696541,-22.19644082133315) rectangle (41.76413067908503,39.5997794785621);

%%%%%%%%%%%%%%%%%%%%%% PROSTE PIONOWE %%%%%%%%%%%%%%%%%%%%%%%%%%%%%%%%%%%%%%%%%%%%%%

%%% R_a
\draw [line width=1pt] (-13.5,30)-- (-13.5,-40);
\draw (-13.5,34) node[anchor=north] {$R_a$};
%%% R_b
\draw [line width=1pt] (0,30)-- (0,-40);
\draw (0,34) node[anchor=north] {$R_b$};
%%% R_c
\draw [line width=1pt] (15,30)-- (15,-40);
\draw (15,34) node[anchor=north] {$R_c$};

%%%%%%%%%%%%%%%%%%%%%% PROSTE UKO君E %%%%%%%%%%%%%%%%%%%%%%%%%%%%%%%%%%%%%%%%%%%%%%

%%% R_d
\draw [line width=1pt] (26,38)-- (16.7,-37);
\draw (26,40) node[anchor=north east] {$R_d$};

%%%%%%%%%%%%%%%%%%%%%%%%%%%%%%%  PUNKTY %%%%%%%%%%%%%%%%%%%%%%%%%%%%%%%%%%%%

\draw [fill=uuuuuu] (-13.5,20) circle (8pt);

\draw [fill=uuuuuu] (-13.5,10) circle (8pt);

\draw [fill=uuuuuu] (-13.5,0) circle (8pt);

\draw [fill=uuuuuu] (-13.5,-10) circle (8pt);

\draw [fill=uuuuuu] (0,20) circle (8pt);

\draw [fill=uuuuuu] (0,10) circle (8pt);

\draw [fill=uuuuuu] (0,0) circle (8pt);

\draw [fill=uuuuuu] (0,-10) circle (8pt);

\draw [fill=uuuuuu] (15,20) circle (8pt);
\draw[color=uuuuuu] (15,20) node[anchor=south east] {$c_1$};

\draw [fill=uuuuuu] (15,10) circle (8pt);
\draw[color=uuuuuu] (15,10) node[anchor=south east] {$c_2$};

\draw [fill=uuuuuu] (15,0) circle (8pt);
\draw[color=uuuuuu] (15,0) node[anchor=south east] {$c_3$};

\draw [fill=uuuuuu] (15,-10) circle (8pt);
\draw[color=uuuuuu] (15,-10) node[anchor=south west] {$c_4$};

\draw [fill=uuuuuu] (24.555276497257866,26.275434462072584) circle (8pt);

\draw [fill=uuuuuu] (23.18572916907802,15.376924033599247) circle (8pt);

\draw [fill=uuuuuu] (21.816181840898172,4.478413605125911) circle (8pt);

\draw [fill=uuuuuu] (19.044287161367107,-17.87962152071977) circle (8pt);

\end{tikzpicture}
    \caption{A (4,4) half grid with numbered points on line $R_c$}
    \label{fig: (4,4) half grid Rc}
\end{figure}
Since $Z_{abc}$ is a grid on $Q_{abc}$, there are four lines $r_1,r_2,r_3,r_4$ in the ruling complementary to that determined by the lines $R_a, R_b$ and $R_c$, each of them passing through the point $c_i$ with the same index. The intersections of each $r_i$ with $R_a, R_b$ determine a labelling for the points of $Z$ on the lines $R_a$ and $R_b$, see Figure \ref{fig: (4,4) half grid Rabc}.
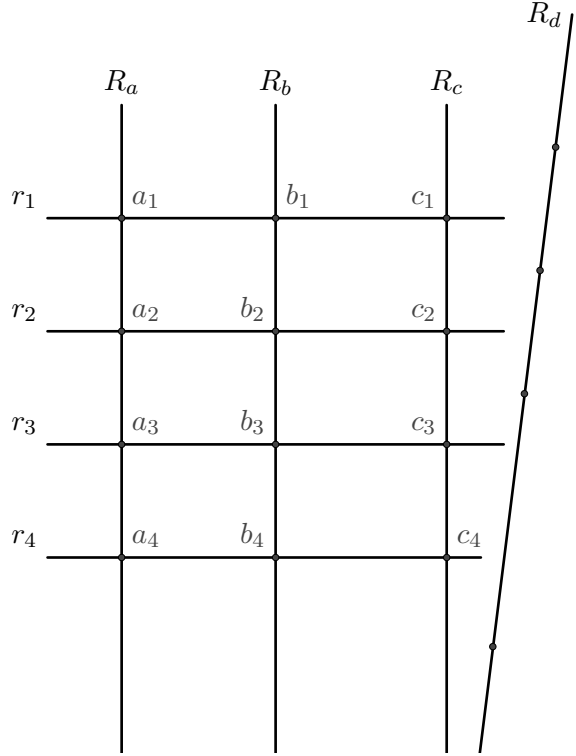
\begin{figure}[h!]
    \centering
\begin{tikzpicture}[line cap=round,line join=round,>=triangle 45,x=1cm,y=1cm,scale=0.15]
\clip(-31.8109674696541,-27.19644082133315) rectangle (41.76413067908503,39.5997794785621);

%%%%%%%%%%%%%%%%%%%%%%%%%%%%%%% PROSTE POZIOME %%%%%%%%%%%%%%%%%%%%%%%%%%%%%%%%%%%%%
%%% r_1
\draw [line width=1pt] (-20,20)-- (20,20);
\draw (-20,20) node[anchor=south east] {$r_1$};
%%% r_2
\draw [line width=1pt] (-20,10)-- (20,10);
\draw (-20,10) node[anchor=south east] {$r_2$};
%%% r_3
\draw [line width=1pt] (-20,0)-- (20,0);
\draw (-20,0) node[anchor=south east] {$r_3$};
%%% r_4
\draw [line width=1pt] (-20,-10)-- (18,-10);
\draw (-20,-10) node[anchor=south east] {$r_4$};

%%%%%%%%%%%%%%%%%%%%%% PROSTE PIONOWE %%%%%%%%%%%%%%%%%%%%%%%%%%%%%%%%%%%%%%%%%%%%%%

%%% R_a
\draw [line width=1pt] (-13.5,30)-- (-13.5,-40);
\draw (-13.5,34) node[anchor=north] {$R_a$};
%%% R_b
\draw [line width=1pt] (0,30)-- (0,-40);
\draw (0,34) node[anchor=north] {$R_b$};
%%% R_c
\draw [line width=1pt] (15,30)-- (15,-40);
\draw (15,34) node[anchor=north] {$R_c$};

%%%%%%%%%%%%%%%%%%%%%% PROSTE UKO君E %%%%%%%%%%%%%%%%%%%%%%%%%%%%%%%%%%%%%%%%%%%%%%

%%% R_d
\draw [line width=1pt] (26,38)-- (16.7,-37);
\draw (26,40) node[anchor=north east] {$R_d$};

%%%%%%%%%%%%%%%%%%%%%%%%%%%%%%%  PUNKTY %%%%%%%%%%%%%%%%%%%%%%%%%%%%%%%%%%%%

\draw [fill=uuuuuu] (-13.5,20) circle (8pt);
\draw[color=uuuuuu] (-13.5,20) node[anchor=south west] {$a_1$};

\draw [fill=uuuuuu] (-13.5,10) circle (8pt);
\draw[color=uuuuuu] (-13.5,10) node[anchor=south west] {$a_2$};

\draw [fill=uuuuuu] (-13.5,0) circle (8pt);
\draw[color=uuuuuu] (-13.5,0) node[anchor=south west] {$a_3$};

\draw [fill=uuuuuu] (-13.5,-10) circle (8pt);
\draw[color=uuuuuu] (-13.5,-10) node[anchor=south west] {$a_4$};

\draw [fill=uuuuuu] (0,20) circle (8pt);
\draw[color=uuuuuu] (0,20) node[anchor=south west] {$b_1$};

\draw [fill=uuuuuu] (0,10) circle (8pt);
\draw[color=uuuuuu] (0,10) node[anchor=south east] {$b_2$};

\draw [fill=uuuuuu] (0,0) circle (8pt);
\draw[color=uuuuuu] (0,0) node[anchor=south east] {$b_3$};

\draw [fill=uuuuuu] (0,-10) circle (8pt);
\draw[color=uuuuuu] (0,-10) node[anchor=south east] {$b_4$};

\draw [fill=uuuuuu] (15,20) circle (8pt);
\draw[color=uuuuuu] (15,20) node[anchor=south east] {$c_1$};

\draw [fill=uuuuuu] (15,10) circle (8pt);
\draw[color=uuuuuu] (15,10) node[anchor=south east] {$c_2$};

\draw [fill=uuuuuu] (15,0) circle (8pt);
\draw[color=uuuuuu] (15,0) node[anchor=south east] {$c_3$};

\draw [fill=uuuuuu] (15,-10) circle (8pt);
\draw[color=uuuuuu] (15,-10) node[anchor=south west] {$c_4$};

\draw [fill=uuuuuu] (24.555276497257866,26.275434462072584) circle (8pt);

\draw [fill=uuuuuu] (23.18572916907802,15.376924033599247) circle (8pt);

\draw [fill=uuuuuu] (21.816181840898172,4.478413605125911) circle (8pt);

\draw [fill=uuuuuu] (19.044287161367107,-17.87962152071977) circle (8pt);

\end{tikzpicture}
    \caption{A $(4,4)$ half grid with numbered points on lines $R_a, R_b$ and $R_c$}
    \label{fig: (4,4) half grid Rabc}
\end{figure}

Since $r_1,\ldots,r_4$ are contained in a quadric we have
\begin{equation}\label{eq: first 3 c-r}
   j(a_1,a_2;a_3,a_4)=j(b_1,b_2;b_3,b_4)=j(c_1,c_2;c_3,c_4)
\end{equation}
by Lemma \ref{lem: quadrics and (2,4) grids}.
This is a considerable constraint and we want to combine it with similar conditions determined by the remaining three quadrics.

\begin{figure}[h!]
    \centering
\begin{tikzpicture}[line cap=round,line join=round,>=triangle 45,x=1cm,y=1cm,scale=0.15]
\clip(-31.8109674696541,-27.19644082133315) rectangle (41.76413067908503,39.5997794785621);

%%%%%%%%%%%%%%%%%%%%%%%%%%%%%%% PROSTE POZIOME %%%%%%%%%%%%%%%%%%%%%%%%%%%%%%%%%%%%%
%%% r_1
\draw [line width=1pt] (-20,20)-- (20,20);
\draw (-20,20) node[anchor=south east] {$r_1$};
%%% r_2
\draw [line width=1pt] (-20,10)-- (20,10);
\draw (-20,10) node[anchor=south east] {$r_2$};
%%% r_3
\draw [line width=1pt] (-20,0)-- (20,0);
\draw (-20,0) node[anchor=south east] {$r_3$};
%%% r_4
\draw [line width=1pt] (-20,-10)-- (18,-10);
\draw (-20,-10) node[anchor=south east] {$r_4$};

%%%%%%%%%%%%%%%%%%%%%% PROSTE PIONOWE %%%%%%%%%%%%%%%%%%%%%%%%%%%%%%%%%%%%%%%%%%%%%%

%%% R_a
\draw [line width=1pt] (-13.5,30)-- (-13.5,-40);
\draw (-13.5,34) node[anchor=north] {$R_a$};
%%% R_b
\draw [line width=1pt] (0,30)-- (0,-40);
\draw (0,34) node[anchor=north] {$R_b$};
%%% R_c
\draw [line width=1pt] (15,30)-- (15,-40);
\draw (15,34) node[anchor=north] {$R_c$};

%%%%%%%%%%%%%%%%%%%%%% PROSTE UKO君E %%%%%%%%%%%%%%%%%%%%%%%%%%%%%%%%%%%%%%%%%%%%%%

%%% R_d
\draw [line width=1pt] (26,38)-- (16.7,-37);
\draw (26,40) node[anchor=north east] {$R_d$};

%%% L_1
\draw [line width=1pt] (-4.735588173233172,6.893152013767983)-- (30.14589050360053,30);
\draw (30.661044078509807,31.358790721144786) node[anchor=north west] {$L_1$};

%%% L_2
\draw [line width=1pt] (-4.6284183563126176,-3.0365377485794993)-- (30.14589050360053,20);
\draw (30.92688242552327,21.12401436112651) node[anchor=north west] {$L_2$};

%%% L_3
\draw [line width=1pt] (-4.396365690770353,-12.884296394334472)-- (30.14589050360053,10);
\draw (30.661044078509807,11.420914695135158) node[anchor=north west] {$L_3$};

%%% L_4
\draw [line width=1pt] (-5.149187653820004,30.134586903116933)-- (29.50818252796565,-38.70757259495852);
\draw (28.8,-35) node[anchor=north east] {$L_4$};

%%%%%%%%%%%%%%%%%%%%%%%%%%%%%%%  PUNKTY %%%%%%%%%%%%%%%%%%%%%%%%%%%%%%%%%%%%

\draw [fill=uuuuuu] (-13.5,20) circle (8pt);
\draw[color=uuuuuu] (-13.5,20) node[anchor=south west] {$a_1$};

\draw [fill=uuuuuu] (-13.5,10) circle (8pt);
\draw[color=uuuuuu] (-13.5,10) node[anchor=south west] {$a_2$};

\draw [fill=uuuuuu] (-13.5,0) circle (8pt);
\draw[color=uuuuuu] (-13.5,0) node[anchor=south west] {$a_3$};

\draw [fill=uuuuuu] (-13.5,-10) circle (8pt);
\draw[color=uuuuuu] (-13.5,-10) node[anchor=south west] {$a_4$};

\draw [fill=uuuuuu] (0,20) circle (8pt);
\draw[color=uuuuuu] (0,20) node[anchor=south west] {$b_1$};

\draw [fill=uuuuuu] (0,10) circle (8pt);
\draw[color=uuuuuu] (0,10) node[anchor=south east] {$b_2$};

\draw [fill=uuuuuu] (0,0) circle (8pt);
\draw[color=uuuuuu] (0,0) node[anchor=south east] {$b_3$};

\draw [fill=uuuuuu] (0,-10) circle (8pt);
\draw[color=uuuuuu] (0,-10) node[anchor=south east] {$b_4$};

\draw [fill=uuuuuu] (15,20) circle (8pt);
\draw[color=uuuuuu] (15,20) node[anchor=south east] {$c_1$};

\draw [fill=uuuuuu] (15,10) circle (8pt);
\draw[color=uuuuuu] (15,10) node[anchor=south east] {$c_2$};

\draw [fill=uuuuuu] (15,0) circle (8pt);
\draw[color=uuuuuu] (15,0) node[anchor=south east] {$c_3$};

\draw [fill=uuuuuu] (15,-10) circle (8pt);
\draw[color=uuuuuu] (15,-10) node[anchor=south west] {$c_4$};

\draw [fill=uuuuuu] (24.555276497257866,26.275434462072584) circle (8pt);
\draw[color=uuuuuu] (24.555276497257866,26.275434462072584) node[anchor=north west] {$d_1$};

\draw [fill=uuuuuu] (23.18572916907802,15.376924033599247) circle (8pt);
\draw[color=uuuuuu] (23.18572916907802,15.376924033599247) node[anchor=north west] {$d_2$};

\draw [fill=uuuuuu] (21.816181840898172,4.478413605125911) circle (8pt);
\draw[color=uuuuuu] (21.816181840898172,4.478413605125911) node[anchor=north west] {$d_3$};

\draw [fill=uuuuuu] (19.044287161367107,-17.87962152071977) circle (8pt);
\draw[color=uuuuuu] (19.244287161367107,-17.87962152071977) node[anchor=west] {$d_4$};

\end{tikzpicture}
    \caption{Grid lines on $Q_{abc}$ and $Q_{bcd}$}
    \label{fig: two grids}
\end{figure}
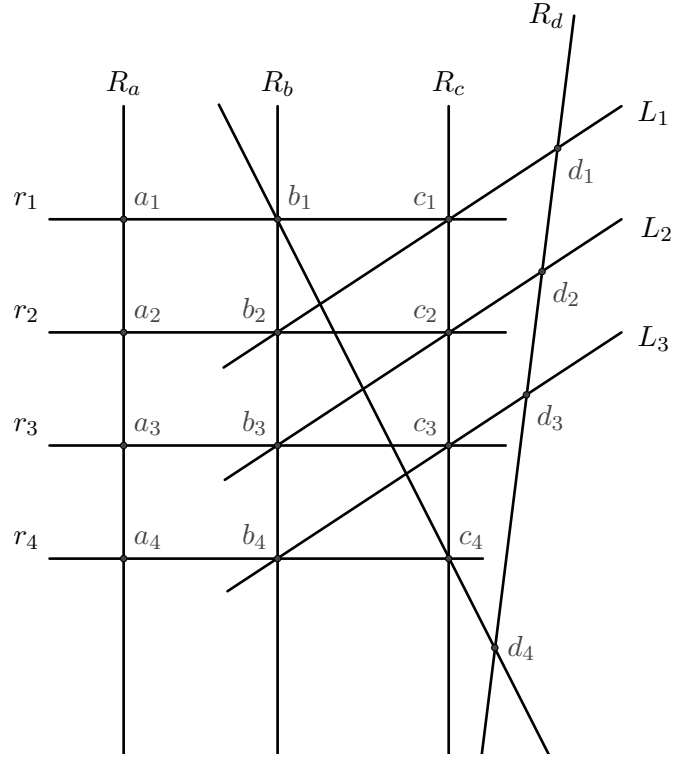

The set $Z_{bcd}$ is a grid in $Q_{bcd}$. So there are four lines in the ruling of $Q_{bcd}$ complementary to the ruling containing $R_c$ cutting out the points of $Z$ on the union $R_b\cup R_c\cup R_d$. We call them $L_1,\ldots,L_4$, where the numbering is determined by the numbering of points on $R_c$, see Figure \ref{fig: two grids}.
This numbering determines also numbering of points on the line $R_d$. Thus all points in $Z$ are now labeled.
For $i=1,\ldots,4$ we denote by $b_{\beta(i)}$ the point of intersection $R_b\cap L_i$ (in Figure \ref{fig: two grids} we took as an example $\beta(1)=2$, $\beta(2)=3$ and so on). In any case we have
\begin{equation}\label{eq: last 3 c-r}
j(b_{\beta(1)},b_{\beta(2)};b_{\beta(3)},b_{\beta(4)})=j(c_1,c_2;c_3,c_4)=j(d_1,d_2;d_3,d_4)
\end{equation}
by Lemma \ref{lem: quadrics and (2,4) grids}.
Now, the key point here is to see that the permutation $\beta$ can be assumed not to be an involution, which forces one of cases: \eqref{eq: harmonic} or \eqref{eq: anharmonic}. We will first exclude the possibility that $\beta$ is the identity.
\begin{proposition}\label{pro: beta not id}
   The permutation $\beta$ is not the identity.
\end{proposition}
\begin{proof}
    If $\beta(i)=i$, then $L_i=r_i$, so that $d_i\in Q_{abc}$. If this happens for all $i$, we obtain that $Z$ is contained in a quadric, hence it is a grid. A contradiction.
\end{proof}
   Since $\beta$ preserves the cross-ratio of four points on $R_b$, there exists a projective transformation $\varphi_\beta$ of $R_b$, which restricts to $\beta$ on $Z\cap R_b$. We show in the next Lemma that this projectivity has exactly two fixed points.
\begin{lemma}[Transversals to $R_a,\ldots,R_d$]\label{lem: transversals}
    There are two distinct transversals $S, S'$ to $R_a,\ldots, R_d$. Moreover, the intersection points of these transversals with $R_b$ are the fixed points of $\varphi_\beta$.
\end{lemma}
\begin{proof}
    The intersection $Q_{abc}\cap Q_{bcd}$ contains two skew lines $R_b$ and $R_c$, so it must contain either one line $S$ (counted with multiplicity $2$) or two lines $S, S'$ from the complementary rulings on both quadrics. Hence these lines must be transversals of $R_a,\ldots,R_d$.

    Let us denote their intersection points:
    $$a_s=S\cap R_a,\; b_s=S\cap R_b,\; c_s=S\cap R_c,\; d_s=S\cap R_d$$
    and similarly
    $$a_{s'}=S'\cap R_a,\; b_{s'}=S'\cap R_b,\; c_{s'}=S'\cap R_c,\; d_{s'}=S'\cap R_d,$$
    see Figure \ref{fig: complete}.
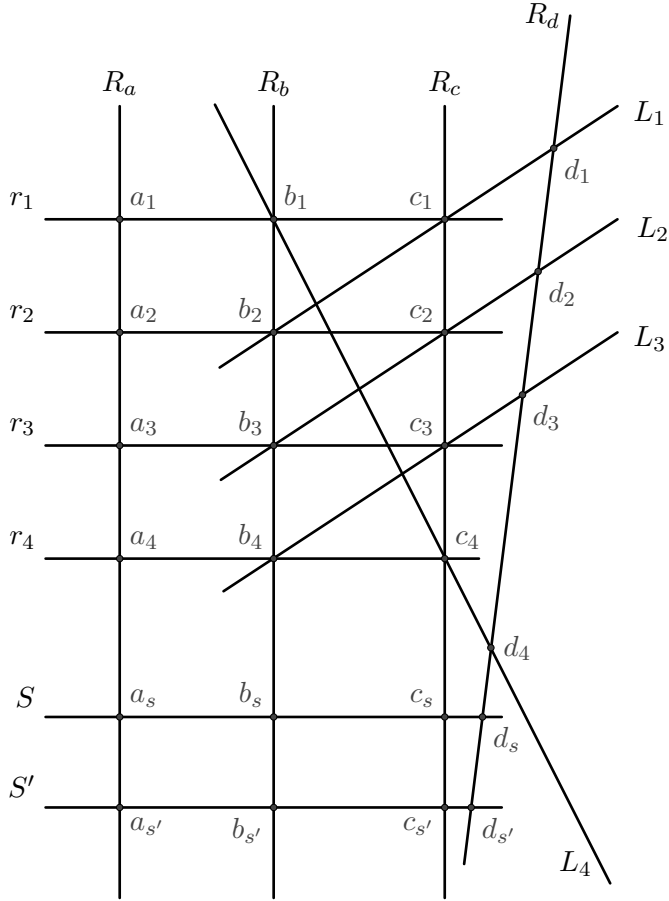
\begin{figure}[h!]
    \centering
\begin{tikzpicture}[line cap=round,line join=round,>=triangle 45,x=1cm,y=1cm,scale=0.15]
\clip(-31.8109674696541,-47.19644082133315) rectangle (61.76413067908503,39.5997794785621);

%%%%%%%%%%%%%%%%%%%%%%%%%%%%%%% PROSTE POZIOME %%%%%%%%%%%%%%%%%%%%%%%%%%%%%%%%%%%%%
%%% r_1
\draw [line width=1pt] (-20,20)-- (20,20);
\draw (-20,20) node[anchor=south east] {$r_1$};
%%% r_2
\draw [line width=1pt] (-20,10)-- (20,10);
\draw (-20,10) node[anchor=south east] {$r_2$};
%%% r_3
\draw [line width=1pt] (-20,0)-- (20,0);
\draw (-20,0) node[anchor=south east] {$r_3$};
%%% r_4
\draw [line width=1pt] (-20,-10)-- (18,-10);
\draw (-20,-10) node[anchor=south east] {$r_4$};
%%% S
\draw [line width=1pt] (-20,-24)-- (20,-24);
\draw (-20,-24) node[anchor=south east] {$S$};
%%% S'
\draw [line width=1pt] (-20,-32)-- (20,-32);
\draw (-20,-32) node[anchor=south east] {$S'$};

%%%%%%%%%%%%%%%%%%%%%% PROSTE PIONOWE %%%%%%%%%%%%%%%%%%%%%%%%%%%%%%%%%%%%%%%%%%%%%%

%%% R_a
\draw [line width=1pt] (-13.5,30)-- (-13.5,-40);
\draw (-13.5,34) node[anchor=north] {$R_a$};
%%% R_b
\draw [line width=1pt] (0,30)-- (0,-40);
\draw (0,34) node[anchor=north] {$R_b$};
%%% R_c
\draw [line width=1pt] (15,30)-- (15,-40);
\draw (15,34) node[anchor=north] {$R_c$};

%%%%%%%%%%%%%%%%%%%%%% PROSTE UKO君E %%%%%%%%%%%%%%%%%%%%%%%%%%%%%%%%%%%%%%%%%%%%%%

%%% R_d
\draw [line width=1pt] (26,38)-- (16.7,-37);
\draw (26,40) node[anchor=north east] {$R_d$};

%%% L_1
\draw [line width=1pt] (-4.735588173233172,6.893152013767983)-- (30.14589050360053,30);
\draw (30.661044078509807,31.358790721144786) node[anchor=north west] {$L_1$};

%%% L_2
\draw [line width=1pt] (-4.6284183563126176,-3.0365377485794993)-- (30.14589050360053,20);
\draw (30.92688242552327,21.12401436112651) node[anchor=north west] {$L_2$};

%%% L_3
\draw [line width=1pt] (-4.396365690770353,-12.884296394334472)-- (30.14589050360053,10);
\draw (30.661044078509807,11.420914695135158) node[anchor=north west] {$L_3$};

%%% L_4
\draw [line width=1pt] (-5.149187653820004,30.134586903116933)-- (29.50818252796565,-38.70757259495852);
\draw (28.8,-35) node[anchor=north east] {$L_4$};

%%%%%%%%%%%%%%%%%%%%%%%%%%%%%%%  PUNKTY %%%%%%%%%%%%%%%%%%%%%%%%%%%%%%%%%%%%

\draw [fill=uuuuuu] (-13.5,20) circle (8pt);
\draw[color=uuuuuu] (-13.5,20) node[anchor=south west] {$a_1$};

\draw [fill=uuuuuu] (-13.5,10) circle (8pt);
\draw[color=uuuuuu] (-13.5,10) node[anchor=south west] {$a_2$};

\draw [fill=uuuuuu] (-13.5,0) circle (8pt);
\draw[color=uuuuuu] (-13.5,0) node[anchor=south west] {$a_3$};

\draw [fill=uuuuuu] (-13.5,-10) circle (8pt);
\draw[color=uuuuuu] (-13.5,-10) node[anchor=south west] {$a_4$};

\draw [fill=uuuuuu] (-13.5,-24) circle (8pt);
\draw[color=uuuuuu] (-13.5,-24) node[anchor=south west] {$a_s$};

\draw [fill=uuuuuu] (-13.5,-32) circle (8pt);
\draw[color=uuuuuu] (-13.5,-32) node[anchor=north west] {$a_{s'}$};

\draw [fill=uuuuuu] (0,20) circle (8pt);
\draw[color=uuuuuu] (0,20) node[anchor=south west] {$b_1$};

\draw [fill=uuuuuu] (0,10) circle (8pt);
\draw[color=uuuuuu] (0,10) node[anchor=south east] {$b_2$};

\draw [fill=uuuuuu] (0,0) circle (8pt);
\draw[color=uuuuuu] (0,0) node[anchor=south east] {$b_3$};

\draw [fill=uuuuuu] (0,-10) circle (8pt);
\draw[color=uuuuuu] (0,-10) node[anchor=south east] {$b_4$};

\draw [fill=uuuuuu] (0,-24) circle (8pt);
\draw[color=uuuuuu] (0,-24) node[anchor=south east] {$b_s$};

\draw [fill=uuuuuu] (0,-32) circle (8pt);
\draw[color=uuuuuu] (0,-32) node[anchor=north east] {$b_{s'}$};

\draw [fill=uuuuuu] (15,20) circle (8pt);
\draw[color=uuuuuu] (15,20) node[anchor=south east] {$c_1$};

\draw [fill=uuuuuu] (15,10) circle (8pt);
\draw[color=uuuuuu] (15,10) node[anchor=south east] {$c_2$};

\draw [fill=uuuuuu] (15,0) circle (8pt);
\draw[color=uuuuuu] (15,0) node[anchor=south east] {$c_3$};

\draw [fill=uuuuuu] (15,-10) circle (8pt);
\draw[color=uuuuuu] (15,-10) node[anchor=south west] {$c_4$};

\draw [fill=uuuuuu] (15,-24) circle (8pt);
\draw[color=uuuuuu] (15,-24) node[anchor=south east] {$c_s$};

\draw [fill=uuuuuu] (15,-32) circle (8pt);
\draw[color=uuuuuu] (15,-32) node[anchor=north east] {$c_{s'}$};

\draw [fill=uuuuuu] (24.555276497257866,26.275434462072584) circle (8pt);
\draw[color=uuuuuu] (24.555276497257866,26.275434462072584) node[anchor=north west] {$d_1$};

\draw [fill=uuuuuu] (23.18572916907802,15.376924033599247) circle (8pt);
\draw[color=uuuuuu] (23.18572916907802,15.376924033599247) node[anchor=north west] {$d_2$};

\draw [fill=uuuuuu] (21.816181840898172,4.478413605125911) circle (8pt);
\draw[color=uuuuuu] (21.816181840898172,4.478413605125911) node[anchor=north west] {$d_3$};

\draw [fill=uuuuuu] (19.044287161367107,-17.87962152071977) circle (8pt);
\draw[color=uuuuuu] (19.244287161367107,-17.87962152071977) node[anchor=west] {$d_4$};

\draw [fill=uuuuuu] (18.309426261993095,-24) circle (8pt);
\draw[color=uuuuuu] (18.349426261993095,-24) node[anchor=north west] {$d_s$};

\draw [fill=uuuuuu] (17.320717721577902,-32) circle (8pt);
\draw[color=uuuuuu] (17.360717721577902,-32) node[anchor=north west] {$d_{s'}$};

\end{tikzpicture}
    \caption{All lines and transversals}
    \label{fig: complete}
\end{figure}
Now we want to exclude the possibility that $S=S'$. To this end we prove the following\\
\textbf{Claim.}\\
The projectivity $\varphi_\beta$ has two fixed points at $b_s$ and $b_{s'}$.\\
\textbf{Proof of the Claim.}\\
By Lemma \ref{lem: projectivities on two lines} there exists a projective transformation $\Phi$ of $\P^3$, which restricts to $\varphi_\beta$ on $R_b$ and to the identity on $R_c$. This projectivity maps lines $r_i$ joining $b_i$ and $c_i$ to lines joining $\Phi(b_i)$ and $\Phi(c_i)$ for $i=1,\ldots,4$. But $\Phi(b_i)=b_{\beta(i)}$ and $\Phi(c_i)=c_i$, so that $\Phi(r_i)=L_i$ for all $i=1,\ldots,4$. It follows that
$$\Phi(Q_{abc})=Q_{bcd}.$$
Moreover, since $\Phi$ leaves $R_b$ and $R_c$ invariant by Lemma \ref{lem: projectivities on two lines} it must leave also the union $S\cup S'$ invariant, because
$$Q_{abc}\cap Q_{bcd}=R_b\cup R_c\cup S \cup S'.$$
It must be in fact $\Phi(S)=S$ because $\Phi$ restricts to the identity on $R_c$ so it cannot swap the points $c_s$ and $c_{s'}$ (the claim $\Phi(S)=S$ remains valid also if $S=S'$).

Since $\varphi_\beta$ is not the identity, it has at most two fixed points. It has also at least two fixed points by Lemma \ref{lem: projectivity of P1 with 1 FP} because it is a transformation of finite order, equal to the order of $\beta$, so that all its orbits are finite. Suppose now that $S=S'$. Then there exists another point $b_s\neq P\in R_b$ fixed by $\varphi_\beta$. Let $r_P$ be the line in the same ruling as $r_1$ on $Q_{abc}$ passing through $P$. This line meets $R_c$ in some point $P_c$. Since $\Phi(P_c)=P_c$ (as $\Phi$ restricted to $R_c$ is the identity) the line $r_P$ is invariant under $\Phi$. But then it belongs to $Q_{abc}\cap Q_{bcd}$, so it must be the line $S'$ distinct from $S$.

The two fixed points of $\varphi_\beta$ are then $b_s$ and $b_{s'}$. This ends the proof of the Claim and also of the Lemma.
\end{proof}
Considering the grid $Z_{abd}$ on the quadric $Q_{abd}$ we obtain another permutation $\beta'$ acting on points $b_1,\ldots,b_4$ determined by numbering of points on the line $R_d$ in a way that the following triples of points are collinear:
\begin{equation}\label{eq: beta'}
a_{\alpha(1)},b_{\beta'(1)},d_1,\;\;
a_{\alpha(2)},b_{\beta'(2)},d_2,\;\;
a_{\alpha(3)},b_{\beta'(3)},d_3,\;\;
a_{\alpha(4)},b_{\beta'(4)},d_4,
\end{equation}
where $\alpha$ is some permutation of points on $R_a$.
\begin{lemma}\label{lem: betas differ}
    The permutations $\beta$ and $\beta'$ do not coincide and at least one of them is not an involution.
\end{lemma}
\begin{proof}
    By definition of $\beta$ the following triples of points are collinear:
\begin{equation}\label{eq: beta}
b_{\beta(1)},c_1,d_1,\;\;
b_{\beta(2)},c_2,d_2,\;\;
b_{\beta(3)},c_3,d_3,\;\;
b_{\beta(4)},c_4,d_4.
\end{equation}
Combining \eqref{eq: beta} with \eqref{eq: beta'} we conclude that if $\beta=\beta'$, then the lines $L_i$ meet the line $R_a$ at points of $Z$. But this implies that $Z$ is a grid, a contradiction.

Now, replacing $R_c$ by $R_d$ in the proof of Lemma \ref{lem: transversals} we conclude that $b_s$ and $b_{s'}$ are fixed points of $\beta'$. So $\beta$ and $\beta'$ induce projectivities $\varphi$ and $\varphi'$ on $R_b$ which have the same pair points as their fixed points. By Lemma \ref{lem: involutions with 2 fixed points} at most one of these projectivities can be an involution.
\end{proof}
Since all permutations in \eqref{eq: standard c-r} are involutions, we conclude that the points in $Z_b$ are either harmonic or anharmonic. Then \eqref{eq: first 3 c-r} and \eqref{eq: last 3 c-r} imply that points in all $Z_x$ are simultaneously either harmonic or anharmonic for $x\in\left\{a,b,c,d\right\}$. We study both cases in more detail in the next section.
%*****************************************************************************
\section{Proof of Theorem \ref{thm: main}}
To fix notation we assume in this section that the permutation $\beta$ introduced before Proposition \ref{pro: beta not id} is not an involution. This can be done since relabelling the lines $R_a, R_b, R_c$ and $R_d$ according to the following rule: $(a,b,c,d)\to (d,b,a,c)$ exchanges the role of $\beta$ and $\beta'$ from Lemma \ref{lem: betas differ}.
\subsection{The anharmonic case}\label{ssec:anharmonic}
In this part we prove Theorem \ref{thm: main} A). We see in \eqref{eq: anharmonic} that all non-involutions appearing there have a fixed point. Renumbering the points, if necessary, we may assume that $b_4$ is the fixed point of $\beta$. This implies that $r_4$ is one of the two lines transversal to $R_a,\ldots,R_d$. Let this line be $S$, so that we have $r_4=S$.

Moreover, renumbering the remaining points, if necessary, we may assume that $\beta$ is the $3$-cycle $(2, 3, 1, 4)$.
Now we can begin with fixing coordinates.
Suppose to begin with that
$$R_b:\; x=z=0\;\mbox{ and }\; R_c:\; x-y=z-w=0.$$
Let also
$$b_1=(0:1:0:0),\; b_2=(0:0:0:1),\; b_3=(0:1:0:1)$$
and
$$c_1=(1:1:0:0),\; c_2=(0:0:1:1),\; c_3=(1:1:1:1).$$
Then necessarily
$$b_4=(0:1:0:\varepsilon)\;\mbox{ and }\; c_4=(1:1:\varepsilon:\varepsilon),$$
where $\varepsilon$ is a primitive root of unity of order $6$, so that it satisfies the equation $\varepsilon^2-\varepsilon+1=0$.
Then it must be
$$r_1:\; z=w=0,\; r_2:\; x=y=0,\;
r_3:\; x-z=y-w=0,\;\mbox{ and }\; r_4=S:\; \varepsilon x-z=\varepsilon y-w=0.$$
With the given fixed $\beta$ we can reproduce also equations of the $L_i$ lines:
$$L_1:\; x-y=z=0,\; L_2:\; x=y+z-w=0,\;
L_3:\; x-z=z-w=0\;\mbox{ and }\; L_4=r_4=S.$$
It remains to construct the lines $R_a$ and $R_d$.

To this end note that for $i\in\left\{1,2,3\right\}$ the line $M_i$ joining $c_i$ to one of points $a_1,a_2,a_3$ on $R_a$ must also meet $R_d$ in one of the points $d_1,d_2,d_3$.
The following Lemma restricts considerably possible collineations.
\begin{lemma}\label{lem: not on M_i}
    It is neither $a_{\beta(i)}\in M_i$ nor $a_i\in M_i$  for all  $i=1,\ldots,4$.
\end{lemma}
\begin{proof}
    Assume to the contrary that $a_{\beta(i)}\in M_i$ for some $i$, see Figure \ref{fig: Lemma not on M_i}.
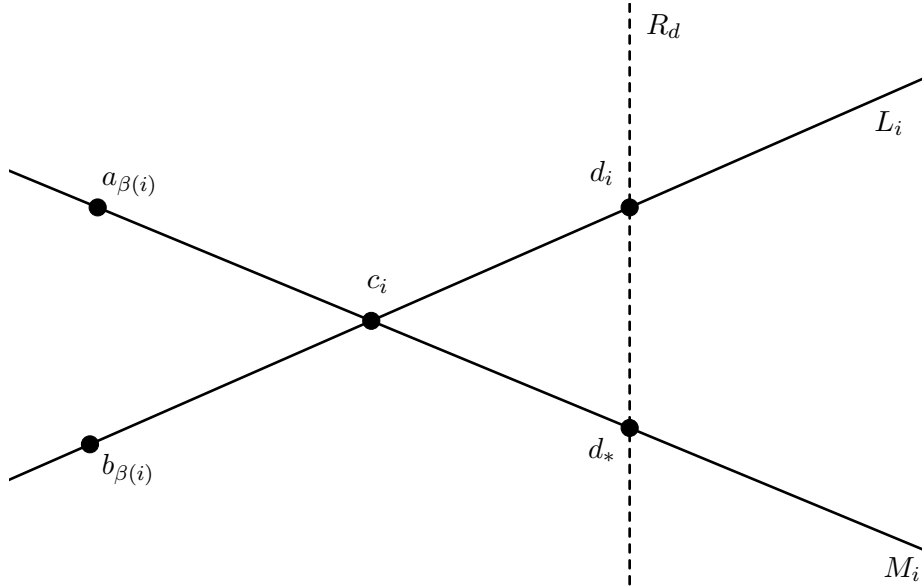
\begin{figure}[h!]
    \centering
\begin{tikzpicture}[line cap=round,line join=round,>=triangle 45,x=1cm,y=1cm, scale=0.8]
\clip(-6,-3.7) rectangle (9,6);
\draw [line width=1pt,domain=-10.28447894019779:16.130991083800332] plot(\x,{(--6.318846651866136-3.6553253955411593*\x)/8.745614383889103});
\draw [line width=1pt,domain=-10.28447894019779:16.130991083800332] plot(\x,{(--6.774347685197417--3.926766407753858*\x)/8.872284119479266});
\draw [line width=1pt,dashed] (4.2,-7.04981222237246) -- (4.2,7.927260876139567);
\draw (4.3,6) node[anchor=north west] {%\large
$R_d$};
\draw (8.056828548163171,4.38858470311343) node[anchor=north west] {%\large
$L_i$};
\draw (8.206350076600897,-3) node[anchor=north west] {%\large
$M_i$};

\draw [fill=black] (-4.545614383889103,2.622405395541159) circle (4pt);
\draw[color=black] (-4.01703487318314,3) node {%\large
$a_{\beta(i)}$};
\draw [fill=black] (4.2,-1.03292) circle (4pt);
\draw[color=black] (3.745624478208761,-1.4) node {%\large
$d_*$};
\draw [fill=black] (-4.672284119479266,-1.3043564077538583) circle (4pt);
\draw[color=black] (-4.041955127922761,-1.7171966891351526) node {%\large
$b_{\beta(i)}$};
\draw [fill=black] (4.2,2.62241) circle (4pt);
\draw[color=black] (3.745624478208761,3.192412475611638) node {%\large
$d_i$};
\draw [fill=black] (-0.047672358678294516,0.7424412225594632) circle (4pt);
\draw[color=black] (0.057426776744872066,1.3476756627966735) node {%\large
$c_i$};
\end{tikzpicture}
    \caption{Ilustration for Lemma \ref{lem: not on M_i}}
    \label{fig: Lemma not on M_i}
\end{figure}
    Since the lines $M_i$ and $L_i$ intersect in $c_i$ (and are not equal, as otherwise $L_i$ would be one of the secants $S,S'$ contradicting the assumption that $\beta$ has no fixed point), they span a plane, call it $\pi_i$. Since $L_i$ and $M_i$ intersect $R_d$ in two distinct points, this line is also contained in $\pi_i$. Similarly, since $b_{\beta(i)}\in L_i$ by definition and $a_{\beta(i)}\in M_i$ by assumption, the line $r_{\beta(i)}$ is contained in $\pi_i$. But then the lines $r_{\beta(i)}$ and $R_d$ intersect, so that $r_{\beta(i)}$ is either $S$ or $S'$, a contradiction.

    It cannot be that $a_i\in M_i$, because then it would be either $M_i=r_i$ and we would have $M_i=S'$ contradicting $\beta$ not being an involution.
\end{proof}
\begin{corollary}\label{cor: point on M_i}
    It is $a_{\beta^2(i)}\in M_i$ for all  $i=1,\ldots,4$.
\end{corollary}
Moreover, replacing $M_i$ by $N_i$ in Lemma \ref{lem: not on M_i} we conclude that neither $a_{\beta^{-1}(i)}\in N_i$ nor $a_i\in N_i$.
\begin{corollary}\label{cor: point on N_i}
    It is $a_{\beta(i)}\in N_i$ for all  $i=1,\ldots,4$.
\end{corollary}

Since we are in the position to pick the point $a_1$ on the $r_1$ line arbitrarily (but distinct from $b_1$ and $c_1$), we pick $a_1=(1:0:0:0)$, this forces $R_a:\; y=w=0$ and
$$a_2=(0:0:1:0),\;
a_3=(1:0:1:0)\mbox{ and } a_4=(1:0:\eps:0).$$
Now, as we have specific coordinates for all points on lines $R_a, R_b, R_c$ we can determine equations of lines $M_i$ and $N_j$ and check their intersections. We obtain
$$d_1=(1:0:1:1),\; d_2=(0:1:-1:0), d_3=(1:1:0:1).$$
Then $R_d:\; y+z-w=x-w=0$ and $d_4=(\eps:\eps-1:1:\eps)$. This concludes proof of Theorem \ref{thm: main} A).

\subsection{The harmonic case}\label{ssec:harmonic}
In this part we prove Theorem \ref{thm: main} B). Right away we fix equations of lines $R_a, R_b, R_c$ and the coordinates of points $a_1,a_2,a_3$, $b_1,b_2,b_3$ and $c_1,c_2,c_3$ as in Section \ref{ssec:anharmonic}. Then necessarily
$$a_4=(1:0:-1:0),\;
b_4=(0:1:0:-1)\;\mbox{ and }\;
c_4=(1:1:-1:-1),$$
since the quadruples of points on lines $R_a,R_b$ and $R_c$ must be harmonic.

The only permutations in \eqref{eq: harmonic} of order greater than $2$ are the $4$-cycles $(3,4, 2, 1)$ and $(4, 3, 1, 2)$. Renumbering the points if necessary, we may assume that $\beta=(3, 4, 2, 1)$. Then the lines $L_i$ are determined as follows:
$$L_1:\; z=x-y+w=0,\;
L_2:\; x=y-z+w=0,\;
L_3:\; x-z=y-z=0\; \mbox{ and }
L_4:\; x+w=z-w=0$$
and
$$c_i,d_i,b_{\beta(i)}\in L_i\;\mbox{ for } i=1,\ldots,4.$$
Let, as before, $M_i$ be the line through $c_i$ containing a configuration point on $R_a$ and another one on $R_d$. And let $N_i$ be the line through $b_i$ containing a configuration point on $R_a$ and on $R_d$. Since Lemma \ref{lem: not on M_i} remains valid also in the situation considered here only the following collinearities are possible:
$$
M_1:\; c_1 \mbox{ and } \left\{ a_2 \mbox{ or } a_4\right\},\;
M_2:\; c_2 \mbox{ and } \left\{a_1 \mbox{ or } a_3\right\},$$$$
M_3:\; c_3 \mbox{ and } \left\{ a_1 \mbox{ or } a_4\right\},\;
M_4:\; c_4 \mbox{ and } \left\{ a_2 \mbox{ or } a_3\right\}.
$$
The choice of points on $M_1$ and the possibilities listed above determine the points on the remaining lines, so that there are only two possibilities:
$$M_1:\; c_1-a_2\; \mbox{ forces }
M_4:\; c_4-a_3 \mbox{ and }
M_2:\; c_2-a_1 \mbox{ and }
M_3:\; c_3-a_4,$$
whereas
$$M_1:\; c_1-c_4\; \mbox{ forces }
M_3:\; c_3-a_1 \mbox{ and }
M_2:\; c_2-a_3 \mbox{ and }
M_4:\; c_4-a_2.$$
By the same token for the lines $N_i$ we obtain the following possibilities
$$
N_1:\; b_1 \mbox{ and } \left\{ a_2 \mbox{ or } a_3\right\},\;
N_2:\; b_2 \mbox{ and } \left\{a_1 \mbox{ or } a_4\right\},$$$$
N_3:\; b_3 \mbox{ and } \left\{ a_2 \mbox{ or } a_4\right\},\;
N_4:\; b_4 \mbox{ and } \left\{ a_1 \mbox{ or } a_3\right\}.
$$
which provides two cases for the configuration:
$$N_1:\; b_1-a_2\; \mbox{ forces }
N_3:\; b_3-a_4 \mbox{ and }
N_2:\; b_2-a_1 \mbox{ and }
N_4:\; b_4-a_3,$$
whereas
$$N_1:\; b_1-a_3\; \mbox{ forces }
N_4:\; b_4-a_1 \mbox{ and }
N_2:\; b_2-a_4 \mbox{ and }
N_3:\; b_3-a_2.$$
We need now to determine which lines $M_i$ and $N_j$ intersect in points being potential configuration points on the $R_d$ line. So these points cannot be the $a_k$ points. Incidences between all possible lines $M_i$ and $N_j$ lines are summarized in Table \ref{tab: incidences}. \begin{table}[h!]
    \centering
    \begin{tabular}{c|c|c|c|c||c|c|c|c|}
    %\hline
    & $b_1a_2$ & $b_2a_1$ & $b_3a_4$ & $b_4a_3$ & $b_1a_3$ & $b_2a_3$ & $b_3a_2$ & $b_4a_1$\\
    \hline
$c_1a_2$ & $a_2$ & $\emptyset$ & $\emptyset$ & $\emptyset$ & $1:1:1:0$ & $\emptyset$ & $a_2$ & $\emptyset$\\
\hline
$c_2a_1$ & $\emptyset$ & $a_1$ & $\emptyset$ & $\emptyset$ & $\emptyset$ & $-1:0:1:1$ & $\emptyset$ & $a_1$ \\
\hline
$c_3a_4$ & $\emptyset$ & $\emptyset$& $a_4$ & $\emptyset$ & $\emptyset$ & $a_4$ & 0:1:2:1 & $\emptyset$\\
\hline
$c_4a_3$ & $\emptyset$ & $\emptyset$ & $\emptyset$ & $a_3$ & $a_3$ & $\emptyset$ & $\emptyset$ & 2:1:0:-1\\
\hline
\hline
$c_1a_4$ & $0:1:1:0$ & $\emptyset$ & $a_4$ & $\emptyset$ & $1:2:1:0$ & $a_4$ & $\emptyset$ & $\emptyset$ \\
\hline
$c_2a_3$ & $\emptyset$ & $1:0:0:-1$ & $\emptyset$ & $a_3$ & $a_3$ & $-1:0:1:2$ & $\emptyset$ & $\emptyset$ \\
\hline
$c_3a_1$ & $\emptyset$ & $a_1$ & -1:1:1:1 & $\emptyset$ & $\emptyset$ & $\emptyset$ & 0:1:1:1 & $a_1$\\
\hline
$c_4a_2$ & $a_2$ & $\emptyset$ & $\emptyset$ & 1:1:1:-1 & $\emptyset$ & $\emptyset$ & $a_2$ & 1:1:0:-1\\
\hline
    \end{tabular}
    \caption{Incidences of potential lines $M_i$ and $N_j$}
    \label{tab: incidences}
\end{table}
An analysis of Table \ref{tab: incidences} together with matching new points to equations of lines $L_i$ shows that it must be
$$R_d: x-z+2w=y-z+w=0$$
and the configuration points on this line are
$$d_1=(2:1:0:-1), d_2=(0:1:2:1), d_3=(1:1:1:0), d_4=(-1:0:1:1)$$
or
$$R_d: x-z+2w=y-z+w=0$$
and the configuration points are
$$d_1=(1:0:0:-1), d_2=(0:1:1:0), d_3=(1:1:1:-1), d_4=(-1:1:1:1).$$
We leave it as an exercise to a motivated reader to check that both configurations obtained this way are projectively equivalent. This ends the proof of Theorem \ref{thm: main}.

%*****************************************************************************

\paragraph*{Acknowledgement.}
Chiantini and Favacchio are members of the Italian GNSAGA-INDAM.
Farnik was partially supported by National Science Centre, Poland, Grant 2018/28/C/ST1/00339.
Harbourne was partially supported by Simons Foundation grant \#524858. Migliore was partially supported by Simons Foundation grant \#839618.
Research of Szemberg and Szpond was partially supported by National Science Centre, Poland, Opus Grant 2019/35/B/ST1/00723.

We are grateful to the organizers and to the INdAM for hosting  the Workshop on The Strong and Weak Lefschetz Properties held in Cortona in the period September 12-16, 2022 and providing stimulating working conditions.

%*****************************************************************************
%\bibliographystyle{abbrv}
%\bibliography{CP}

%\end{document}

%\begin{thebibliography}{99}\footnotesize\itemsep=0cm

%***************************************************************************** % Addresses

\bigskip
\small

\noindent
   Luca Chiantini,\\
   Dipartimento di Ingegneria dell'Informazione e Scienze Matemati\-che, Universit\`a di Siena, Italy

\nopagebreak
\noindent
   \textit{E-mail address:} \texttt{luca.chiantini@unisi.it}\\

\bigskip
\noindent
   {\L}ucja Farnik,\\
   Department of Mathematics, University of the National Education Commission, Krakow,
   Podchor\c a\.zych 2,
   PL-30-084 Krak\'ow, Poland

\nopagebreak
\noindent
   \textit{E-mail address:} \texttt{lucja.farnik@gmail.com}\\

\bigskip
\noindent
   Giuseppe Favacchio,\\
   Dipartimento di Ingegneria, Universit\`a degli studi di Palermo,
Viale delle Scienze,  90128 Palermo, Italy

\nopagebreak
\noindent
   \textit{E-mail address:} \texttt{giuseppe.favacchio@unipa.it}\\

\bigskip
\noindent
   Brian Harbourne,\\
   Department of Mathematics,
University of Nebraska,
Lincoln, NE 68588-0130 USA

\nopagebreak
\noindent
   \textit{E-mail address:} \texttt{brianharbourne@unl.edu}\\

\bigskip
\noindent
   Juan Migliore,\\
   Department of Mathematics,
University of Notre Dame,
Notre Dame, IN 46556 USA

\nopagebreak
\noindent
   \textit{E-mail address:} \texttt{migliore.1@nd.edu}\\

\bigskip

\noindent
   Tomasz Szemberg,\\
   Department of Mathematics, University of the National Education Commission, Krakow,
   Podchor\c a\.zych 2,
   PL-30-084 Krak\'ow, Poland

\nopagebreak
\noindent
   \textit{E-mail address:} \texttt{tomasz.szemberg@gmail.com}\\

\bigskip
\noindent
   Justyna Szpond,\\
   Department of Mathematics, University of the National Education Commission, Krakow,
   Podchor\c a\.zych 2,
   PL-30-084 Krak\'ow, Poland.

\nopagebreak
\noindent
   \textit{E-mail address:} \texttt{szpond@gmail.com}\\

%*****************************************************************************

\end{document}